\documentclass{amsart}
\usepackage{amsmath,amssymb,epsfig,graphics,psfrag}

\theoremstyle{plain}
\newtheorem{theorem}{Theorem}
\newtheorem{lemma}[theorem]{Lemma}
\newtheorem{corollary}[theorem]{Corollary}
\newtheorem{proposition}[theorem]{Proposition}

\theoremstyle{definition}
\newtheorem{definition}[theorem]{Definition}

\newcommand{\DONOTUSE}[1]{}
\newcommand{\diam}{{\mathrm{diam}}}
\newcommand{\vol}{{\mathrm{vol}}}
\newcommand{\Rr}{{\mathbb R}} 
\newcommand{\0}{\emptyset}

\renewcommand{\vec}{}

\newcommand{\includefigure}[3]{{
  \begin{center}
  \resizebox{#1}{#2}{\includegraphics{{#3}}}
  \end{center}}}

\allowdisplaybreaks

\begin{document}

\title[Random Geometric Graph Diameter]{Random Geometric Graph
Diameter in the Unit Ball}

\author{Robert B.\ Ellis}
\address{Department of Applied Mathematics\\
Illinois Institute of Technology\\
Chicago, IL 60616\\
USA} \email{rellis@math.iit.edu}

\author{Jeremy L.\ Martin}
\address{Department of Mathematics\\
University of Kansas\\
Lawrence, KS 66045-7523}
\email{jmartin@math.ku.edu}

\author{Catherine Yan}
\address{Department of Mathematics\\
Texas A\&M University\\
College Station, TX 77843-3368\\
USA}
\email{cyan@math.tamu.edu}

\begin{abstract}
The {\em unit ball random geometric graph} $G=G^d_p(\lambda,n)$
has as its vertices $n$ points distributed independently and
uniformly in the unit ball in ${\mathbb R}^d$, with two vertices
adjacent if and only if their $\ell_p$-distance is at most
$\lambda$. Like its cousin the Erd\H{o}s-R\'enyi random graph, $G$
has a \emph{connectivity threshold}: an asymptotic value for
$\lambda$ in terms of $n$, above which $G$ is connected and below
which $G$ is disconnected. In the connected zone, we determine
upper and lower bounds for the graph diameter of $G$.
Specifically, almost always, $\diam_p(\mathbf{B})(1-o(1))/\lambda
\leq \diam(G) \leq \diam_p(\mathbf{B})(1+O((\ln \ln n/\ln
n)^{1/d}))/\lambda$, where $\diam_p(\mathbf{B})$ is the
$\ell_p$-diameter of the unit ball $\mathbf{B}$. We employ a
combination of methods from probabilistic combinatorics and
stochastic geometry.
\end{abstract}

\maketitle

\section{Introduction}

A \textit{random geometric graph} consists of a set of vertices
distributed randomly over some metric space $X$, with two vertices
joined by an edge if the distance between them is sufficiently
small.  This construction presents a natural alternative to the
classical Erd\H{o}s-R\'enyi random graph model, in which the
presence of each edge is an independent event (see, e.g.,
\cite{B01}). The study of random geometric graphs is a relatively
new area; the monograph \cite{P03} by M.~Penrose is the current
authority. In addition to their theoretical interest, random
geometric graphs have many applications, including wireless communication
networks; see, e.g.,
\cite{chen-jia:2001,stojmenovic-others:2002,wu-li:2002}.

In this article, we study the \textit{unit ball random geometric
graph} $G=G^d_p(\lambda,n)$, defined as follows. Let $d$ and $n$
be positive integers, $\mathbf{B}$ the Euclidean unit ball in
$\Rr^d$ centered at the origin, $\lambda$ a positive real number,
and $p \in [1,\infty]$ (that is, either $p \in [1,\infty)$ or
$p=\infty$). Let $V_n$ be a set of $n$ points in $\mathbf{B}$,
distributed independently and uniformly with respect to Lebesgue
measure on ${\mathbb R}^d$. Then $G$ is the graph with vertex set
$V_n$, where two vertices $\vec{x}=(x_1,\dots,x_d)$ and
$\vec{y}=(y_1,\dots,y_d)$ are adjacent if and only if
$\|\vec{x}-\vec{y}\|_p\leq\lambda$. (Thus the larger $\lambda$ is,
the more edges $G$ has.) Here $\|\cdot\|_p$ is the
\textit{$\ell_p$-metric} defined by
  $$\|\vec{x}-\vec{y}\|_p = \begin{cases}
    \big(\sum_{i=1}^d |x_i-y_i|^p\big)^{1/p} & \quad\text{for } p\in[1,\infty),\\
    \max\{|x_i-y_i| ~:~ 1 \leq i \leq d\}    & \quad\text{for }
    p=\infty,
  \end{cases}$$
where the case $p=2$ gives the standard Euclidean metric on
$\Rr^d$.

When $d=1$, $G$ is known as a \emph{random interval graph}. (Note
that the value of $p$ is immaterial when $d=1$.)  Random interval
graphs have been studied extensively in the literature; the
asymptotic distributions for the number of isolated vertices and
the number of connected components were determined precisely by
E.~Godehardt and J.~Jaworski \cite{GJ96}. The \emph{random
Euclidean unit disk graph} $G^2_2(\lambda,n)$ was studied by
X.~Jia and the first and third authors \cite{EJY04}.

In the present article, we focus on the case $d\geq 2$ and $p\in
[1,\infty]$, but also comment along the way on the special case
$d=1$. We are
interested in the asymptotic behavior of the connectivity and
graph diameter of $G$ as $n\to\infty$ and $\lambda\to 0$. In fact,
$G$ has a \emph{connectivity threshold}: roughly speaking, an
expression for $\lambda$ as a function of $n$, above which $G$ is
connected and below which $G$ is disconnected. (This
behavior is ubiquitous in the theory of the Erd\H{o}s-R\'enyi
random graph model; cf.\ \cite{B01}.)

We are interested primarily in the combinatorial graph diameter of
$G$, $\diam(G)$, above the connectivity threshold.
Our results include

\begin{itemize}
  \item a lower bound for $\diam(G)$
  (Proposition~\ref{prop:diameter-lower-bound} of \S\ref{sec:lowerBound});
  \item an ``absolute'' upper bound $\diam(G) < K/\lambda$, where $K$ is a
    constant depending only on $d$ (Theorem~\ref{thm:diameter}
    of \S\ref{sec:absoluteUpperBound}); and
  \item an asymptotically tight upper bound within a factor of the form $(1+o(1))$
    of the lower bound, the proof of which builds on the
    absolute upper bound (Theorem~\ref{thm:diameterBound} of
    \S\ref{sec:upperBound}).
\end{itemize}

\section{Definitions and Notation}

As mentioned above, the main object of our study is the random geometric
graph $G=G_p^d(\lambda,n)$, where $d\geq 1$ is the dimension of the ambient unit ball
${\mathbf B}$,
$p\in[1,\infty]$ describes the metric, $\lambda>0$ is the $\ell_p$-distance
determining adjacency, and $n$ is the number of vertices.
We will generally avoid repeating the constraints on the parameters.

The graph distance $d_G(\vec{x},\vec{y})$ between two vertices
$\vec{x},\vec{y}\in V_n$ is defined to be the length of the
shortest path between $\vec{x}$ and $\vec{y}$ in $G$, or $\infty$
if there is no such path. The graph diameter of $G$ is defined to
be $\diam(G):=\max\{d_G(\vec{x},\vec{y}) ~:~ \vec{x},\vec{y}\in
V_n\}$.  This graph-theoretic quantity is not to be confused with
the $\ell_p$-diameter of a set $X\subseteq\mathbb{R}^d$, defined
as $\diam_p(X):=\sup\{\|\vec{x}-\vec{y}\|_p ~:~ \vec{x},\vec{y}\in
X\}$.  The $\ell_p$-ball of radius $r$ centered at
$\vec{x}\in\mathbb{R}^d$ is defined as
$$ B_p^d(\vec{x},r) := \{\vec{y}\in\mathbb{R}^d ~:~
    \|\vec{x}-\vec{y}\|_p\leq r\},
$$
while the $\ell_p$-ball of radius $r$ around a set $X\subseteq
\mathbb{R}^d$ is $B_p^d(X,r):=\cup_{x\in X}B_p^d(x,r)$. The origin
of $\mathbb{R}^d$ is $O:=(0,\ldots,0)$; when the center of a ball
is not explicitly given, we define $B_p^d(r):=B_p^d(O,r)$. Thus
$\mathbf{B}=B_2^d(1)$. The $\ell_p$-diameter of $\mathbf{B}$ is
$$
\diam_p(\mathbf{B}) ~=~
\begin{cases}
2d^{1/p-1/2} &   \mbox{ when } 1\leq p\leq 2,\\
2 & \mbox{ when } 2\leq p\leq \infty.
\end{cases}
$$
The distance $d_p(X,Y)$ between two sets $X,Y\subseteq
\mathbb{R}^d$ is defined as $\inf\{\|x-y\|_p~:~x\in X,y\in Y\}$.
The \emph{boundary} $\partial X$ of $X$ is its closure minus its
interior (in the usual topology on $\mathbb{R}^d$), and its
\emph{volume} $\vol(X)$ is its Lebesgue measure.

We will make frequent use of the quantity
  \begin{equation} \label{eqn:alpha}
  \alpha_p^d ~:=~ \frac{\vol(B_p^d(r))}{\vol(B_2^d(r))} ~=~
  \frac{\Gamma\left(\frac{p+1}{p}\right)^d \cdot \Gamma\left(\frac{2+d}{2}\right)}
       {\Gamma\left(\frac{3}{2}\right)^d   \cdot \Gamma\left(\frac{p+d}{p}\right)}\,,
  \end{equation}
where $\Gamma$ is the usual gamma function (see, e.g.,
\cite{WhitWat}). The calculation of $\alpha_p^d$, along with the
proofs of several other useful facts about $\ell_p$-geometry, may
be found in the Appendix at the end of the article.

We will say that the random graph $G$ has a property $P$
\emph{almost always}, or \emph{a.a.}, if
  $$
  \lim_{n\rightarrow \infty} \Pr\left[G \text{ has property $P$}\right] = 1.
  $$
By the notation $f(n)=o(g(n))$ and $f(n)=O(g(n))$, we mean,
respectively,  $\lim_{n\rightarrow\infty}f(n)/g(n)=0$ and
$\limsup_{n\rightarrow\infty}f(n)/g(n)\leq c$, for some absolute
nonnegative constant $c$.

\section{Connectivity thresholds}

In order for $G=G_p^d(\lambda,n)$ to have finite diameter, it must
be connected.  Therefore, we seek a \emph{connectivity threshold}
-- a lower bound on $\lambda$ so that $G$ is almost always
connected. When $d=1$, all $\ell_p$-metrics are identical.  For
this case we now quote parts of Theorems 10 and 12 of \cite{GJ96},
to which we refer the reader for their precise determination of
the asymptotic Poisson distributions of the number of isolated
vertices and the number of connected components.

\begin{theorem}[Godehardt, Jaworski]\label{thm:d1Thresh}
Let $\lambda_1(n)=\frac{1}{n}(\log n + c + o(1))$ and
$\lambda_2(n)=\frac{2}{n}(\log n + c + o(1))$, where $c$ is a
constant. Then
  \begin{eqnarray}
  \lim_{n\rightarrow\infty}\mathrm{Pr}[G_p^1(\lambda_1(n),n)
    \mbox{has an isolated vertex}\,]& = & e^{-e^{-c}}\,,\nonumber\\
  \lim_{n\rightarrow\infty}\mathrm{Pr}[G_p^1(\lambda_2(n),n)
    \mbox{is connected}\,]& = & e^{-e^{-c}}\,.\nonumber
\end{eqnarray}
\end{theorem}

In particular, by replacing $c$ in Theorem \ref{thm:d1Thresh} with
a nonnegative sequence $\gamma(n)\rightarrow\infty$, almost always
$G_p^1(\lambda_1(n),n)$ has no isolated vertices and
$G_p^1(\lambda_2(n),n)$ is connected. The case $d=1$ is
exceptional in that the thresholds for having isolated vertices
and for connectivity are separated.

For $d\geq 2$, we will use the fact that the connectivity
threshold coincides with the threshold for the disappearance of
isolated vertices, which follows from two theorems of M.~Penrose.
First we compute the threshold for isolated vertices, which is
easier to calculate.

\begin{proposition}\label{prop:UBIsolatedVertices}
Let $d\geq 2$, let $p\in[1,\infty]$, and let $\alpha=\alpha^d_p$
be the constant of \eqref{eqn:alpha}.  Suppose $\gamma(n)$ is a
nonnegative sequence such that
$\lim_{n\rightarrow\infty}\gamma(n)\rightarrow\infty$, and that
$$
\lambda ~\geq ~
\left(\frac{1}{\alpha n}
    \left(\frac{2(d-1)}{d}\ln n+\frac{2}{d}\ln\ln
    n+\gamma(n)\right)\right)^{1/d}\,.
$$
Then, almost always, $G=G_p^d(\lambda,n)$ has no isolated vertices.
\end{proposition}

\begin{proof}
Let $V_n=\{v_1, v_2, \dots, v_n\}$ be the vertex set of $G$. For
each vertex $v_i$, let $A_i$ be the event that $v_i$ is an isolated
vertex, and let $X_i$ be the indicator of $A_i$; that is, $X_i=1$
if $A_i$ occurs and $0$ otherwise.  Set $X=X_1+X_2+\cdots+X_n$.  We
will show that $\mathbb{E}[X]=o(1)$.

By definition, $v_i$ is isolated if and only if there are no
other vertices in $B_p^d(v_i, \lambda) \cap \mathbf{B}$. We
condition $\mathrm{Pr}[A_i]$ on the $\ell_2$-distance from $v_i$
to the origin $O$.
 If $\|v_i-O\|_2\in [0,1-d^{1/2}\lambda)$, then
$B_p^d(v_i,\lambda)\subseteq \mathbf{B}$. Otherwise, if
$\|v_i-O\|_2\in (1-d^{1/2}\lambda,1]$, then the volume of
$B_p^d(v_i, \lambda) \cap \mathbf{B}$ is not less than
$\frac12\vol(B_p^d(v_i, \lambda))(1+O(\lambda))$. Hence
\begin{eqnarray*} \label{pra_i}
\Pr[A_i] & \leq &
    (1-d^{1/2}\lambda)^d(1-\alpha\lambda^d)^{n-1}
    +  \nonumber \\
& & (1-(1-d^{1/2}\lambda)^d)
\left(1-\frac{\alpha\lambda^d}{2}(1+O(\lambda))\right)^{n-1}.
\end{eqnarray*}
Using $1-x = e^{-x}(1+o(1))$ as $x\to0$ and the binomial
expansion, we have
\begin{eqnarray*}
\Pr[A_i]&\le& (1+o(1)) e^{-\alpha n\lambda^d}
+(dd^{1/2}\lambda+O(\lambda^2))(1+o(1)) e^{-\alpha n\lambda^d/2}.
\end{eqnarray*}
The first term is $o(n^{-1})$ for $d\geq 2$. By linearity of
expectation, $\mathbb{E}[X]=n\cdot\mathrm{Pr}[A_i]$, and so
\begin{eqnarray*}
\mathbb{E}[X]&\le& o(1) +
d^{3/2}n\lambda(1+o(1))n^{-1+1/d}(\ln n)^{-1/d}e^{-\gamma(n)/2}\,.
\end{eqnarray*}
The second term is $o(1)$,
and so $\mathrm{Pr}[X>0]\leq\mathbb{E}[X]\leq o(1)$; that is,
almost always, $G$ has no isolated vertices.
\end{proof}

The number of isolated vertices below the threshold
is easy to compute in certain special cases. For
example, if $p\in[1,\infty]$, $d=2$, $\lambda=\sqrt{c\ln n/n}$,
and $0\leq c<\alpha^{-1}$, then a minor modification of
\cite[Theorem 1]{EJY04} yields $X=(1+o(1))n^{1-\alpha c}$
almost always.  In general, to determine the behavior
of $X$ more exactly would require complicated integrals
that describe the volume of $B_p^d(v_i, \lambda) \cap \mathbf{B}$ near
the boundary of $\mathbf{B}$ (cf.\ \cite[Chapter 8]{P03}).
For our purposes, it suffices to concentrate on the values of $\lambda$ for
which $G$ has no isolated vertices.

For $d\geq 2$ and $p\in (1,\infty]$, the connectivity threshold
for the unit-\emph{cube} random geometric graph coincides with the
threshold for lacking isolated vertices.  We quote Penrose's
theorem \cite[Thm.~1.1]{penrose:1999} after some supporting
definitions. Define the \emph{unit cube geometric graph} $H=H^d_p(\lambda,n)$
analogously to $G=G^d_p(\lambda,n)$, except that its vertices are points in
$[0,1]^d$ rather than $\mathbf{B}$.  For any nonnegative
integer $k$, define
  \begin{subequations}
  \begin{align}
  \rho(H; \kappa>k+1) &= \min\{\lambda ~|~ H \text{ has vertex
          connectivity } \kappa\geq k+1\}, \label{eqn:rho1}\\
  \rho(H; \delta>k+1) &= \min\{\lambda ~|~ H \text{ has minimum
          degree } \delta\geq k+1\}. \label{eqn:rho2}
  \end{align}
  \end{subequations}

\begin{theorem}[Penrose]\label{thm:hittingTimes}
Let $p\in (1,\infty]$ and let $k\geq 0$ be an integer.  Then
$$ \lim_{n\rightarrow\infty}\mathrm{Pr}[\,
    \rho(H; \kappa\geq k+1)=
        \rho(H; \delta\geq k+1)\,] \ = \ 1.$$
\end{theorem}

When $k=0$, Theorem~\ref{thm:hittingTimes} asserts that as
$\lambda$ increases (forcing more edges into the graph), almost
always, $H$ becomes connected simultaneously as the last isolated
vertex disappears. In the proof of Theorem \ref{thm:hittingTimes}
in \cite{penrose:1999}, Penrose shows that the limiting
probability distributions for $\rho(H; \kappa\geq k+1)$ and
$\rho(H; \delta\geq k+1)$ are the same. The proof requires only a
series of geometric and probabilistic arguments which hold in the
unit ball as well as in the unit cube (see, in particular,
Sections 2 and 5 of \cite{penrose:1999}), so we have as an
immediate corollary the following.

\begin{corollary}\label{cor:hittingTimes}
Let $d\geq 2$ and $p\in(1,\infty]$, and let $\lambda=\lambda(n)$ be
sufficiently large so that, almost always, $G=G_p^d(\lambda(n),n)$
has no isolated vertices.  Then, almost always, $G$ is connected.
\end{corollary}

We now consider the case that $d\geq 2$ and $p=1$.  Here Theorem 3
does not apply.  However, we can appeal to two general results
about the behavior of a random geometric graph in an
$\ell_p$-metric space whose boundary is a compact
$(d-1)$-submanifold of $\Rr^d$ For such a graph, Theorem 7.2 of
\cite{P03} provides a threshold for the disappearance of isolated
vertices, and Theorem 13.7 provides a threshold for connectivity.
Applying these results to $G$, with the thresholds for $G$ defined
as in \eqref{eqn:rho1} and \eqref{eqn:rho2}, we obtain the
following fact.

\begin{proposition}\label{prop:weakThreshold}
Let $d\geq 2$, $p\in[1,\infty]$, and $G=G_p^d(\lambda(n),n)$.  Let
$\alpha=\alpha_p^d$ be the constant of \eqref{eqn:alpha}, and let
$k\geq 0$ be an integer.  Then, almost always,
  $$
  \lim_{n\rightarrow\infty} \left(\frac{n\alpha}{\log n}\rho(G;\kappa\geq k+1)^d\right) ~=~
  \lim_{n\rightarrow\infty} \left(\frac{n\alpha}{\log n}\rho(G;\delta\geq k+1)^d\right) ~=~
  \frac{2(d-1)}{d}\,.
  $$
\end{proposition}

We now collect the above results to present the connectivity
thresholds that we will use in the rest of the paper.

\begin{theorem}\label{thm:conThresh} Let $G=G_p^d(\lambda,n)$ and let
$\alpha=\alpha_p^d$ be the constant of \eqref{eqn:alpha}.
\renewcommand{\labelenumi}{(\roman{enumi})}
  \begin{enumerate}
  \item Suppose that $d\geq 2$, $p\in (1,\infty]$, $\gamma(n)$ is a
nonnegative sequence such that $\gamma(n)\rightarrow\infty$ and
  $$ 
  \lambda ~\geq ~ \left(\frac{1}{\alpha n}
    \left(\frac{2(d-1)}{d}\ln n+\frac{2}{d}\ln\ln
    n+\gamma(n)\right)\right)^{1/d}.
  $$
Then, almost always, $G$ is connected.

\item Suppose that $d\geq 2$, $p\in [1,\infty]$, and
$\lambda=(c\ln n/n)^{1/d}$ for some constant $c>0$.  Then, almost
always, $G$ is connected if $c>(2(d-1)/(d\alpha))$, and
disconnected if $c<(2(d-1)/(d\alpha))$.

\item Suppose that $d=1$, $p\in [1,\infty]$, and $\lambda= 2(\ln
n+\gamma(n))/n$. Then, almost always, $G$ is connected if
$\gamma(n)\rightarrow\infty$ and disconnected if
$\gamma(n)\rightarrow-\infty$.
\end{enumerate}
\renewcommand{\labelenumi}{(\arabic{enumi})}
\end{theorem}

Assertion~(i) follows from combining Proposition~\ref{prop:UBIsolatedVertices} with
Corollary~\ref{cor:hittingTimes}, and assertion~(ii) is implied by
Proposition~\ref{prop:weakThreshold}. (When $p>1$ and $d\geq 2$, the lower bound in (ii) is
implied by the stronger bound in (i).) Assertion~(iii) is implied by
Theorem~\ref{thm:d1Thresh} and its accompanying remarks.

\section{A lower bound for diameter\label{sec:lowerBound}}

When $G$ is connected, $\mathbf{B}$ will usually contain two
vertices whose $\ell_p$-distance is (asymptotically)
$\diam_p(\mathbf{B})$. Therefore, the diameter of $G$ will almost
always be at least $\diam_p(\mathbf{B})(1-o(1))/\lambda$.  The
precise statement is as follows.

\begin{proposition}[Diameter lower bound]\label{prop:diameter-lower-bound}
Let $d\geq 1$ and $p\in[1,\infty]$, and
suppose that $\lambda=\lambda(n)$ is sufficiently large so that
Theorem \ref{thm:conThresh} guarantees that almost always,
$G=G_p^d(\lambda,n)$ is connected.
If $h=h(n)$ satisfies
  \begin{equation} \label{cap-height-growth}
    \lim_{n\rightarrow\infty} h^{(d+1)/2}n = \infty,
  \end{equation}
then, almost always,
  $$\diam(G) ~\geq~ \frac{1-h}{\lambda}\, \diam_p(\mathbf{B}) ~=~
    \begin{cases}
      2(1-h)d^{1/p-1/2}/\lambda & \quad\text{if } p\leq 2,\\
      2(1-h)           /\lambda & \quad\text{if } p\geq 2.
  \end{cases}$$
\end{proposition}

\begin{proof}
Let $\pm\vec{a}$ be a pair of antipodes of the unit ball
$\mathbf{B}$, chosen as in Figure \ref{fig:antipodes},
and let $\pm C$ be the
spherical cap formed by slicing $\mathbf{B}$ with hyperplanes at
distance $h$ from $\pm\vec{a}$ respectively, perpendicular to the
line joining $\vec{a}$ and $-\vec{a}$. Let $A$ be the event that
at least one of the two caps $\pm C$ contains no vertex of $V_n$.
Then
  $$
  \Pr[A] ~=~ 2\Pr[C\cap V_n=\0]
  ~=~ 2\left(1-\frac{\vol(C)}{\vol(\mathbf{B})}\right)^n
  ~\leq~ 2\,\exp\left(-n\frac{\vol(C)}{\vol(\mathbf{B})}\right).
  $$
On the other hand, $\vol(C)/\vol(B)=O(h^{(d+1)/2})$ by
\eqref{cap-bound} of \S\ref{appendix:spherical-cap-subsection},
which together with the condition \eqref{cap-height-growth} on $h$
implies that $\Pr[A]=o(1)$. That is, $G$ almost always contains a
vertex in each of $C$ and $-C$. The result now follows from the
definition of $G$ and the lower bound
\eqref{cap-distance-formula} on the $\ell_p$-distance
between $C$ and $-C$.
\end{proof}

Note that for all $d\geq 1$, $h$ can be chosen to satisfy both
(\ref{cap-height-growth}) and $\lim_{n\rightarrow\infty}
h/\lambda=0$. Also, if the limit in \eqref{cap-height-growth} is a
nonnegative constant, then $\lim_{n\rightarrow\infty}\Pr[A]>0$;
that is, vertices are not guaranteed in both caps.  For the case
$p=2$, Proposition~\ref{prop:diameter-lower-bound} can be
strengthened by identifying a collection of mutually disjoint
antipodal pairs of caps of height $h$ and showing that, almost
always, both caps in at least one pair contain a vertex. Such a
collection corresponds to an antipodally symmetric spherical code
(see \cite{CS98}).

\section{The absolute upper bound\label{sec:absoluteUpperBound}}

In this section we prove that when $G$ is connected, the graph
distance $d_G(x,y)$ between two vertices $x,y\in V_n$ is at most
$K\|x-y\|_p/\lambda$, where $K>0$ is a constant independent of $n$
and $p$, but dependent on $d$.  As a consequence, $\diam(G)\leq
K\diam_p(\mathbf{B})/\lambda$.  This will not be strong enough to
meet (asymptotically) the lower bound in Proposition
\ref{prop:diameter-lower-bound}, but does guarantee a short path
between any pair of vertices.  This fact will be used repeatedly
in the proof of the tight upper bound in Theorem
\ref{thm:diameterBound} of \S\ref{sec:upperBound}. It is
sufficient to prove the following Theorem \ref{thm:diameter},
since for any two points $x, y \in \mathbb{R}^d$, we have $\|x-y\|_2 \leq
d^{1/2}\|x-y\|_p$.

\begin{theorem} \label{thm:diameter}
Let $d\geq 2$, and suppose that $\lambda=\lambda(n)$ is
sufficiently large so that Theorem \ref{thm:conThresh} guarantees
that almost always, $G=G_p^d(\lambda,n)$ is connected. Then for
any two points $x,y\in V_n$ there exists a constant $K$
independent of $n$ and $p$ such that as $n\rightarrow \infty$,
almost always,
$$  d_G(x,y) ~\leq ~ \frac{K\|x-y\|_2}{\lambda}. $$
\end{theorem}

The proof is based on Proposition \ref{prop:LozDiam} below. For
any two vertices $x,y\in V_n$, let
\[
T_{x,y}(k)=\big[ \mbox{convex closure of } (B_2^d(x, k\lambda)
\cup B_2^d(y, k\lambda)) \big]\cap \mathbf{B}.
\]
Thus $T_{x,y}$ is a ``lozenge''-shaped region.
Let $A_n(k)$ be the event that there exist two vertices $x,y \in
V_n$ such that (i) at least one point is inside $B_2^d(O,
1-(k+\sqrt{d})\lambda)$, and (ii) there is no path of $G$ that
lies in $T_{x,y}(k)$ and connects $x$ and $y$.  The proof of
our next result uses ingredients from \cite[p.~285]{P03},
adapted and extended for our present purposes.

\begin{proposition} \label{prop:LozDiam}
Under the same assumptions as in Theorem \ref{thm:diameter}, there
exists a constant $k_0>0$, such that for all $k>k_0$,
\[ \lim_{n
\rightarrow \infty} \Pr[A_n(k)] =0.
\]
\end{proposition}

\begin{proof}
First, we cover the unit ball $\mathbf{B}$ with $d$-dimensional cubes, each of
side length $\epsilon\lambda$, where $\epsilon=1/(4d)$.  Let $L_d$
be the set of centers of these cubes, and for each $z \in L_d$,
denote the closed cube centered at $z$ by $Q_z$.

Suppose $A_n(k)$ occurs for a pair of vertices $x,y$;
without loss of generality, suppose $y\in B_2^d(O,1-(k+\sqrt{d})\lambda)$.
Abbreviate $T_{x,y}(k)$ by $T_{x,y}$.

{\bf Step 1.} First we construct a connected subset $P\subseteq
T_{x,y}$ such that
\renewcommand{\labelenumi}{(\roman{enumi})}
\begin{enumerate}
\item $B_p^d(P,\lambda/4) \subseteq \mathbf{B}$;
\item $\mathrm{diam}_2 (P) \geq (k-\sqrt{d})\lambda$; and
\item $B_p^d(P,\lambda/4)\cap V_n  = \emptyset$.
\end{enumerate}
\renewcommand{\labelenumi}{(\arabic{enumi})}

Let $V_T=V_n\cap T_{x,y}$ be the set of vertices of $G$ lying in
$T_{x,y}$. Then $B_p^d(V_T,\lambda/2)$ is (topologically)
disconnected with $x$ and $y$ lying in different connected
components. Let $D_x$ be the connected component of
$B_p^d(V_T,\lambda/2)$ containing $x$. Let $S$ be the closure of
the connected component of $T_{x,y}\setminus D_x$  containing $y$.
Let $T$ be the closure of $T_{x,y}\setminus S$, so that $T$
contains $x$.  Then both $S$ and $T$ are connected, and their
union is $T_{x,y}$.  The lozenge $T_{x,y}$ is simply connected, so
it is unicoherent \cite[Lemma~9.1]{P03}; by definition of
unicoherence, since $T_{x,y}$ is the union of closed connected
sets $S,T\subseteq T_{x,y}$, then $P_1:=S\cap T$ is connected.
Since $x\in T$, $y\in S$, and $P_1$ separates $x$ and $y$, any
path in $T_{x,y}$ from $x$ to $y$ must pass through $P_1$. In
particular, $P_1$ intersects the line segment joining $x$ and $y$.
Let $u$ be one of the intersection points.

Next, we show that there is a point $w$ on $\partial T_{x,y}$ such
that $d_p(P_1, w) \leq \lambda/2$, and  derive from this that
$\mathrm{diam}_2 (P_1) \geq (k-\sqrt{d}/2) \lambda$ provided that
the $\ell_2$-distance between $u$ and $w$ is at least $k\lambda$
(see Figure \ref{fig:linel} for an illustration). To achieve this,
we must avoid the case that $w$ lies in the boundary of
$\mathbf{B}$.  To this end, let $C_1 =
\partial\mathbf{B}\cap T_{x,y}$, and let $C_2=\partial
T_{x,y}\setminus C_1$. If $d_p(P_1,C_2)\geq \lambda/2$, then $C_2$
must be a subset either of $S$ or of $T$; without loss of
generality, assume $S$.  Then $y$ is disconnected from $x$ in $G$,
which happens with probability tending to zero by Theorem
\ref{thm:conThresh}. Hence, almost always,
$d_p(P_1,C_2)<\lambda/2$.  It follows that there is a point $w \in
C_2$ such that  $d_p(P_1, w) \leq \lambda/2$. Furthermore,
$B_p^d(P_1,\lambda/2)\cap V_n=\emptyset$, by definition of $P_1$
as the intersection of $S$ and $T$.

\begin{figure}
\begin{center}
 \psfrag{bdB}{$\partial\mathbf{B}$}
 \psfrag{bdT}{$\partial T_{x,y}$}
 \psfrag{P}{$P_1$}
 \centerline{\epsfxsize=1.5in \epsfbox{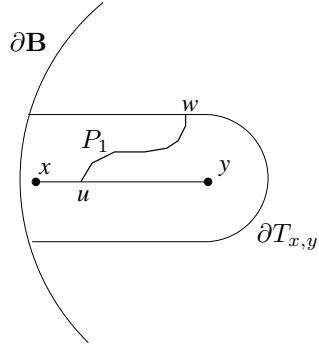}}
 \end{center}\caption{The frontier $P_1$ must intersect the line segment between
 vertices $x,y\in V_n$ at some point $u$, and must also satisfy
 $d_p(P_1,w)\leq \lambda/2$ for some point $w$ on the boundary
 $\partial T_{x,y}\setminus \partial \mathbf{B}$.
\label{fig:linel}}
\end{figure}

As constructed, $P_1$  may be too close to the boundary of
$\mathbf{B}$ so that some cube $Q_z$ intersecting $P_1$ might not
lie entirely inside $\mathbf{B}$. To overcome this, we let $P$ be
obtained from $P_1$ by moving every point toward $O$ by
$\lambda/4$ under the transformation $x\rightarrow
x-(\lambda/4)(x/\|x\|_2)$. Then $P$ is connected, and
$\mathrm{diam}_2(P)\geq \mathrm{diam}_2(P_1)-\lambda/2\geq
(k-\sqrt{d})\lambda$; that is $P$ satisfies conditions (i--iii).

{\textbf{Step 2}}. We now show that when $k$ is large enough, the
probability
  $$\Pr[Q_z\cap V_n=\0 \text{ for every } Q_z \subseteq B_p^d(P,\lambda/4)]$$
tends to zero.  Let $\omega$ be the set of points $z\in L_d$ such
that $Q_z\cap P\neq\emptyset$. Since $P$ is connected, $\omega$ is
a $*$-connected subset of $L_d$; that is, the union of the
corresponding set of cubes is (topologically) connected (see
Figure \ref{fig:lozenge}). For each $z\in \omega$, we have
$Q_z\cap P\neq \emptyset$ and $\epsilon \leq 1/(4d^{1/p})$; hence
$Q_z\subseteq B_p^d(P,\lambda/4)$. By considering the
$\ell_2$-diameter of $P$, we see that $\omega$ contains at least
$4\sqrt{d} (k-\sqrt{d})$ points.  Hence we have a $*$-connected
subset $\omega\subseteq L_d$ with cardinality at least $4\sqrt{d}
(k-\sqrt{d})$ such that $Q_z\cap V_n=\emptyset$ for all
$z\in\omega$. We show that the probability of such an event is
$o(1)$.

\begin{figure}
\begin{center}
 \centerline{\epsfxsize=2.5in \epsfbox{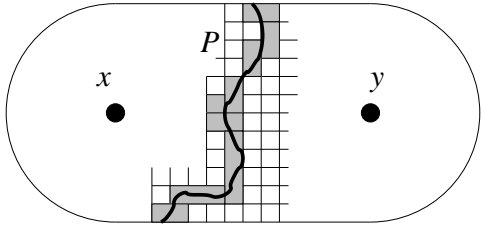}}
\caption{Two vertices $x,y\in V_n$ which are not connected by any
path in $T_{x,y}(k)$, and the ``frontier'' $P$ separating them,
when $d=2$.  The gray squares are the $*$-connected subset
$\omega$ intersecting $P$ of the set of squares covering
$\mathbf{B}$. \label{fig:lozenge}}
 \end{center}
\end{figure}

Let $\mathcal{C}_{d,i}$ denote the collection of $*$-connected
sets of $\omega\subseteq L_d$ of cardinality $i$. It is known that
the number of $*$-connected subsets of $\mathbb{Z}^d$ of
cardinality $i$ containing the origin is at most $2^{3^d i}$ (see,
for example, \cite[Lemma 9.3]{P03}). Since $|L_d|\leq (2/(\epsilon
\lambda))^d$, we have $\mathcal{C}_{d,i} \leq (2/(\epsilon
\lambda))^d 2^{3^d i} \leq 2^d (\epsilon \lambda)^{-d} e^{3^d i}$.
Therefore
\begin{eqnarray}
\Pr[A_n(k) ] & \leq & \sum_{i \geq
4\sqrt{d}\left(k-\sqrt{d})\right) } \ \sum_{\omega \in
\mathcal{C}_{d,i} }
 \Pr[ V_n \cap ( \cup_{z \in \omega} Q_z) =\emptyset ] \nonumber \\
& \leq & \sum_{i \geq 4\sqrt{d}\left(k-\sqrt{d})\right)} 2^d
(\epsilon\lambda)^{-d} \exp(3^d i)
 \left(1-\frac{i}{\mathrm{vol}(\mathbf{B})} (\epsilon\lambda)^d \right)^n
 \nonumber  \\
& \leq & \sum_{i \geq 4\sqrt{d}\left(k-\sqrt{d})\right)} c n
\exp(-i\epsilon^d(d-1)\ln n / \left(d\alpha\,
    \mathrm{vol}(\mathbf{B})\right)) \label{ineq}
\\
& = &  O\left(n^{1-\left(4\sqrt{d}(k-\sqrt{d})(d-1)
\right)/\left(d\alpha\,\mathrm{vol}(\mathbf{B})
(4d)^d\right)}\right), \label{quantity}
\end{eqnarray}
where $c$ is a constant and $\alpha=\alpha_p^d$ is the constant of
\eqref{eqn:alpha}. To justify inequality \eqref{ineq}, when $n$ is
sufficiently large, we have $3^d < \epsilon^d (d-1)\ln n /
(d\alpha\,\mathrm{vol}(\mathbf{B})).$ The order bound in
\eqref{quantity} is immediate by geometric series, and the
resulting quantity is $o(1)$ provided $k>\sqrt{d} +
d\alpha\,\mathrm{vol}(\mathbf{B})(4d)^d/(4\sqrt{d}(d-1))$, which
proves the existence of $k_0$ in the proposition.
\end{proof}

\begin{proof}[Proof of Theorem \ref{thm:diameter}]
Fix $k > k_0$ as in Proposition \ref{prop:LozDiam}. Let $x$ and
$y$ be two vertices in $V_n$ with $\| x-y\|_p > \lambda$. If at
least one of $x,y$ lies in $B_2^d(O, 1-(k+\sqrt{d})\lambda)$,
then, almost always, there is a path of $G$ connecting $x$
and $y$ in $T_{x,y}$. Suppose the shortest path between $x$ and
$y$ in $T_{x,y}$ has length $g$. Then the $\ell_p$-balls of radius
$\lambda/2$ around every other vertex in the path must be pairwise disjoint,
and each must lie inside the convex closure of $B_2^d(x,
(k+\sqrt{d}/2)\lambda) \cup B_2^d(y, (k+\sqrt{d}/2)\lambda)$.  By
comparing the volume of the $\ell_p$-balls of radius $\lambda/2$
to the volume of $T_{x,y}$, we obtain
\begin{eqnarray*}
\left\lfloor \frac{g}{2}\right\rfloor
    \mathrm{vol}\left(B_p^d(x,\lambda/2)\right) & \leq &
    \mathrm{vol}\left(B_2^d(x,(k+\sqrt{d}/2)\lambda)\right)
    \\
& &    + d_2(x,y)\cdot
\mathrm{vol}\left(B_2^{d-1}(x,(k+\sqrt{d}/2)\lambda)\right),
\end{eqnarray*}
which implies that $g\leq K_1+K_2d_2(x,y)/\lambda \leq
(K_1\sqrt{d}+K_2)d_2(x,y)/\lambda$, where $K_1$ and $K_2$ are
constants independent of $n$ and $p$.

If both $x$ and $y$ lie outside $B_2^d(O, 1-(k+\sqrt{d})\lambda)$,
then we can travel from $x$ to an intermediate vertex $x_1$ just
inside $B(O,1-(k+\sqrt{d})\lambda)$ via a path of bounded length,
and then on to $y$. To this end, let
$r=\max\{(\alpha_p^d)^{1/d}\sqrt{d}, \sqrt{d}\}$, and let $E_n(k)$
be the event that there is a vertex $z \in V_n$ such that $z
\not\in B(O,1-(k+\sqrt{d})\lambda)$ and $V_n \cap (
B(O,1-(k+\sqrt{d})\lambda) \cap B(z,(k+\sqrt{d}+2r)\lambda))
=\emptyset$. Then
\[
\Pr[E_n(k)] \leq n\left(1-(1-(k+\sqrt{d})\lambda)^d\right)
(1-(\sqrt{d}\lambda)^d )^n  =o(1).
\]
Applying this observation with $z=x$, we can find a point $x_1 \in
V_n$ inside $B(O,1-(k+\sqrt{d})\lambda) \cap
B(x,(k+\sqrt{d}+2r)\lambda)$. By the preceding argument we can
first travel from $x$ to $x_1$ in $K d_2(x, x_1)/\lambda$ steps,
and then from $x_1$ to $y$ in $K d_2(x_1, y)/ \lambda$ steps. The
total length of the path is no more than $K (d_2(x, y)+2d_2(x,
x_1))/\lambda$. Theorem \ref{thm:diameter} follows from the fact
that $d_2(x, x_1) \leq (k+\sqrt{d} +2r)\lambda$.
\end{proof}

We briefly discuss the case that $d=1$, so that ${\mathbf B}$ is
the interval $[-1,1]\subset\Rr$. Suppose that $\lambda=\lambda(n)$
is sufficiently large so that Theorem \ref{thm:conThresh}
guarantees that, almost always, $G_p^1(\lambda,n)$ is connected.
For any two vertices $x,y$, the shortest path between them clearly
consists of a strictly increasing set of vertices
$x=x_0<x_1<x_2<\cdots<x_{d_G(x,y)}=y$.  Moreover, the balls
$B(x_0,\lambda/2)$, $B(x_2,\lambda/2)$, $B(x_4,\lambda/2), \dots$
must be pairwise disjoint (else some $x_i$ is redundant). Hence
$|x-y|\geq \lceil d_G(x,y)/2 \rceil \lambda$, and from this it is
not hard to deduce that  $d_G(x,y)\leq 2|x-y|/\lambda$.

\section{The asymptotically tight upper bound\label{sec:upperBound}}

In this section, we improve the upper bound in Theorem
\ref{thm:diameter}, reducing the constant $K$ to
$\diam_p(\mathbf{B})$ (asymptotically). Our main result is as
follows:

\begin{theorem}\label{thm:diameterBound}
Let $d\geq 2$ and $p\in[1,\infty]$, and suppose that
$\lambda=\lambda(n)$ is sufficiently large so that Theorem
\ref{thm:conThresh} guarantees that almost always,
$G=G_p^d(\lambda,n)$ is connected.
Then as $n\to\infty$, almost always,
  $$
  \diam(G) \leq \begin{cases}
    \left(2d^{1/p-1/2}+O\left((\ln \ln n/\ln n)^{1/d}\right)\right)/\lambda
      & \quad\text{ when } 1\leq p\leq 2,\\
    \left(2+O\left((\ln \ln n/\ln n)^{1/d}\right)\right)/\lambda
      & \quad\text{ when } 2\leq p\leq \infty.
  \end{cases}
  $$
That is, almost always, $\diam(G) \leq
\diam_p(\mathbf{B})(1+O((\ln \ln n/\ln n)^{1/d}))/\lambda$.
\end{theorem}

\begin{figure}
\includefigure{3.75in}{2in}{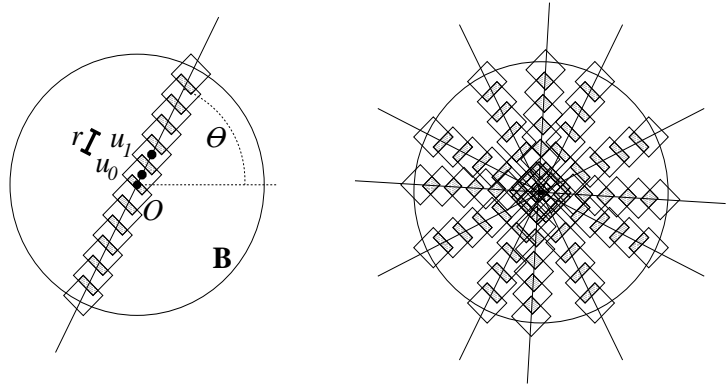}
\caption{{\bf{(a)}} A pin in the unit circle in $\Rr^2$, with
angle $\theta$.  Here $p=1$, so that the $\ell_p$-balls are
diamonds, and  $r$ is the $\ell_2$-distance between consecutive
centers (such as $u_0$ and $u_1$) of $\ell_1$-balls. The petals
are the shaded regions. {\bf{(b)}} A pincushion consisting of
several pins.
\label{fig:pinAndCushion}}
\end{figure}

The proof uses the geometric ingredients of \emph{pins} and
\emph{pincushions}.  A pin consists of a collection of evenly
spaced, overlapping $\ell_p$-balls whose centers lie on a diameter
of the Euclidean unit $d$-ball $\mathbf{B}$. By making suitable
choices for the geometry, we can ensure that each intersection of
consecutive balls contains a vertex in $V_n$, so that the pin
provides a ``highway'' through $G$.  Having done this, we
construct a pincushion so that every point of $\mathbf{B}$ is
reasonably close to an $\ell_p$-ball in one of its constituent
pins.  The following definitions are illustrated in Figure
\ref{fig:pinAndCushion}.

\begin{definition}[Pins]
Fix $d\geq 2$, $p\in[1,\infty]$, $\theta\in(-\pi/2,\pi/2]$,
$\underline{\phi}=(\phi_3,\ldots,\phi_d)\in[0,\pi/2]^d$, and $\lambda,r>0$.  For $m
\in \mathbb{Z}$, put
  $$u_m = u_m(r,\theta) =
    \left(\frac{r}{2}+rm\right)\cdot
    \left( x_1,\ldots,x_d\right)\in \mathbb{R}^d,
  $$
where
  \begin{align*}
    x_1 &= \cos\theta\prod_{i=3}^d\sin\phi_i, &
    x_2 &= \sin\theta\prod_{i=3}^d\sin\phi_i, &
    x_j &= \cos\phi_j\prod_{i=j+1}^d \sin\phi_i \text{ for } 3\leq j\leq d.
  \end{align*}
The corresponding \emph{pin}
$U(d,p,\theta,\underline{\phi},r,\lambda)$ consists of the points
$\{u_m:m\in\mathbb{Z}\}
\cap \mathbf{B}$, together with a collection of $\ell_p$-balls of
radius $\lambda/2$, one centered at each point $u_m$.  Note that
the total number of $\ell_p$-balls is $1+2\lfloor\frac{1}{r}\rfloor$.
\end{definition}

\begin{definition}[Pincushions]
Fix $d\geq 2$, $1\leq p\leq \infty$, and
$\sigma\in\mathbb{Z}^+$.  The corresponding \emph{pincushion}
(with parameters $d,p,\sigma,r,\lambda$)
is the set of $(2\sigma)^{d-1}$ pins
  $$\mathcal{U} := \left\{U(d,p,\theta,\underline{\phi},r,\lambda) ~:~
    \theta,\phi_i \in \left\{ 0, \frac{\pi}{2\sigma}, \frac{2\pi}{2\sigma}, \dots,
    \frac{(2\sigma-1)\pi}{2\sigma} \right\} \right\}.$$
\smallskip
\end{definition}

\begin{definition}[Petals]
Let $U$ be a pin.  A \emph{petal} is the region of intersection of
two overlapping $\ell_p$-balls on $U$ (the shaded regions in
Figure~\ref{fig:pinAndCushion}). A petal is \emph{nonempty} if it
contains a vertex of $V_n$.
\end{definition}

The probability that a petal is nonempty
depends on its volume, which depends in turn on the parameters of the
corresponding pin.  Certainly, we must choose $r$
so that the petal has positive volume: for
example, it suffices to take $r\leq
(\diam_2(\mathbf{B})/\diam_p(\mathbf{B}))\lambda$.
Unfortunately, the volume is difficult to calculate exactly.
Even finding the minimum volume over all angles, that is,
  $$\xi = \xi_p^d(r,\lambda/2) := \inf\{\vol\left(B_p^d(x,\lambda/2)\cap
    B_p^d(y,\lambda/2)\right) ~:~ x,y\in \mathbf{B},\
    \|x-y\|_2=r\}.
  $$
requires integrals that are not easily evaluated
(although for fixed $d$ and $r$, it is certainly true that
$\xi=\Theta(\lambda^d)$).  The easiest way to
find a lower bound for $\xi$ is to inscribe another $\ell_p$-ball
in the petal.

\begin{lemma}\label{lem:volIntBound}
Let $d\geq 1$, $p\in[1,\infty]$, $\lambda>0$ and $0\leq r\leq
(\diam_p(\mathbf{B})/\diam_2(\mathbf{B}))\lambda$. Then, for all
$x,y\in \mathbf{B}$ with $\|x-y\|_2=r$,
  $$\textstyle
  B_2^d\left(\frac{x+y}{2},r'\right)\subseteq
  B_p^d(x,\lambda/2)\cap B_p^d(y,\lambda/2),
  $$
where
  $$
  r' = \begin{cases}
    \frac{\lambda}{2}d^{1/2-1/p}-\frac{r}{2}
      & \quad\mbox{when } 1\leq p\leq 2,\\
    \frac{\lambda}{2}-\frac{r}{2}
      & \quad\mbox{when } 2\leq p \leq \infty.
    \end{cases}
  $$
In particular,
  $$
  \frac{\xi}{\vol(\mathbf{B})} ~\geq~
  \begin{cases}
    \left(\frac{\lambda}{2}d^{1/2-1/p}-\frac{r}{2}\right)^d &
      \mbox{when } 1\leq p\leq 2,\\
    \left(\frac{\lambda}{2}-\frac{r}{2}\right)^d
      & \mbox{when } 2\leq p \leq \infty.
  \end{cases}
  $$
\end{lemma}

\begin{proof}
Let $z\in B^d_2(\frac{x+y}{2},r')$.  By the triangle inequality,
  $$\textstyle
    \|z-x\|_2 ~\leq~ \left\|z-\frac{x+y}{2}\right\|_2 + \left\|\frac{x+y}{2}-x\right\|_2
    ~\leq~ r'+\frac{r}{2},
  $$
which, together with \eqref{lp-diameter}, implies that $\|z-x\|_p \leq \lambda/2$.
That is, $z \in B_p^d(x,\lambda/2)$.  The same argument implies that
$z \in B_p^d(y,\lambda/2)$.  The bound on $\xi$ is then a simple application of
\eqref{lp-ball-volume}.
\end{proof}

\smallskip

For the remainder of this section, we work with
the pincushion $\mathcal{U}$ defined by
  \begin{subequations}
  \begin{equation} \label{eqn:defnSigma}
    \sigma=\sigma(n)=\lfloor(\ln n)^{1/d}\rfloor
  \end{equation}
and
  \begin{equation}\label{eqn:defnR}
    r=r(n)=\begin{cases}
      \lambda d^{1/2-1/p}(1-\rho(n)) & \mbox{ when } 1\leq p\leq 2,\\
      \lambda(1-\rho(n)) & \mbox{ when } 2\leq p\leq \infty ,
    \end{cases}
  \end{equation}
where
  \begin{equation}\label{eqn:defnRho}
    \rho=\rho(n)=
    \begin{cases}
      2d^{1/p-1/2}(\ln \ln n/\ln n)^{1/d}
        & \mbox{ when } 1\leq p\leq 2,\\
    2(\ln \ln n/\ln n)^{1/d}
        & \mbox{ when } 2\leq p\leq \infty.
    \end{cases}
  \end{equation}
  \end{subequations}
Let $\tau_U=\tau_U(n)$ be the number of empty petals along the pin $U\in \mathcal{U}$,
and define
  $$
  \tau = \tau(n) = \max\{\tau_U(n) ~:~ U\in\mathcal{U}\}.
  $$
We first calculate an upper bound on $\tau$.

\begin{lemma}\label{lem:tauBound}
With the assumptions of Theorem \ref{thm:diameterBound}, and the
parameters $\sigma,r,\rho$ as just defined,
almost always,
  $$
  \tau \leq \sigma^{(d-1)/2}\,\frac{2}{r}
  \exp{\left(\frac{-n\xi}{\vol(\mathbf{B})}\right)}.
  $$
\end{lemma}

\begin{proof}
Denote the right-hand side of the desired inequality by $T$.
By linearity of expectation,
\begin{eqnarray}
\mathrm{Pr}[\tau\geq T] & \leq & \mathbb{E}[|\{U:
    \tau_U\geq T\}|] \nonumber \\
& \leq & (2\sigma)^{d-1}\cdot \mathrm{Pr}[\tau_{U^*}\geq T],
    \label{eqn:tauUB}
\end{eqnarray}
where $U^*$ is chosen so as to maximize
$\Pr[\tau_{U}\geq T]$.  Let $X_i$ be
the indicator random variable of the event that the $i$th petal of
$U^*$ is empty, and let $X=\sum_i X_i$. Now $U^*$ contains at most
$2/r$ petals, so by linearity of expectation,
  \begin{equation} \label{eqn:tauExp}
  \mathbb{E}[X] ~\leq~ \frac{2}{r}
    \left(1-\frac{\xi}{\vol(\mathbf{B})} \right)^n
      ~\leq~ \frac{2}{r} \exp{\left(\frac{-n\xi}{\vol(\mathbf{B})}\right)}.
  \end{equation}
Also,
\begin{eqnarray}\label{eqn:tauVar}
\mathrm{var}[X] & \leq & \mathbb{E}[X] +
    \sum_{i\neq j}\mathrm{cov}(X_i,X_j),
\end{eqnarray}
where the covariance $\mathrm{cov}(X_i,X_j)$ for $i\neq j$ is
\begin{eqnarray}
\mathrm{cov}(X_i,X_j) & =& \Pr[X_iX_j=1]-
    \Pr[X_i=1]\cdot\Pr[X_j=1] \nonumber \\
    \left(1-\frac{\xi}{\vol(\mathbf{B})} \right)^{2n}
        \nonumber \\
& \leq & o(1)\cdot \exp{\left(\frac{-2n\xi}{\vol(\mathbf{B})}\right)}.
    \label{eqn:tauCov}
\end{eqnarray}
Combining \eqref{eqn:tauExp}, \eqref{eqn:tauVar} and \eqref{eqn:tauCov} gives
  \begin{equation}\label{eqn:tauVarTwo}
  \mathrm{var}[X] \leq 
  \frac{2}{r}
    \exp\left(\frac{-n\xi}{\vol(\mathbf{B})}\right) + \frac{4}{r^2}\;o(1)
    \cdot \exp{\left(\frac{-2n\xi}{\vol(\mathbf{B})}\right)}.
  \end{equation}
By Chebyshev's inequality (cf.~\cite{Alonspencer:2000}) and the
bounds for $\mathbb{E}[X]$ and $\mathrm{var}[X]$ in
(\ref{eqn:tauExp}) and (\ref{eqn:tauVarTwo}),
  \begin{eqnarray*}
  \Pr\big[X \geq T\big] & \leq &
  \Pr\big[|X-\mathbb{E}[X]|\geq T-\mathbb{E}[X]\big]\\
  &\leq&
    \frac{\mathrm{var}[X]}{(T-\mathbb{E}[X])^2}
  ~=~ o\left(\frac{1}{\sigma^{d-1}}\right) .
  \end{eqnarray*}
Now $\Pr[\tau_U\geq T]=\Pr[X\geq T]$.  Therefore, substituting
this last bound into (\ref{eqn:tauUB}) gives
$\mathrm{Pr}[\tau\geq T] = o(1)$, which implies the desired result.
\end{proof}
\smallskip

We can now prove the main result of this section.

\begin{proof}[Proof of Theorem \ref{thm:diameterBound}]
Let $x,y\in V_n$. We will find vertices $x_1$ and $y_1$ near $x$
and $y$ respectively, belonging to petals of the pincushion
$\mathcal{U}$.  From each of $x_1,y_1$, we walk along the
appropriate pin to points $x_2,y_2$ near the origin and belonging
to petals on the same pin as $x_1$ and $y_1$, respectively. We
will then use Theorem \ref{thm:diameter} to construct a path from
$x_2$ to $y_2$, as well as any ``detours'' needed in case there
are missing edges in the paths along the pins.  Without loss of
generality, we may assume that
  \begin{equation} \label{wolog-assumption}
     \|x\|_2,\:\|y\|_2 ~\geq~ (\tau+3/2)r.
  \end{equation}
The justification for this is deferred until the end of the proof.

Let $R_x$ be the (possibly empty) petal nearest to $x$, and let
$U_x$ be the pin containing $R_x$. When $n$ is sufficiently large,
the distance from $x$ to $R_x$ is at most $d\pi/2\sigma$.  By
definition of $\tau$, there is another vertex $x_1\in V_n$, which
also lies in a petal on $U_x$, but is closer to the origin, so
that $\|x-x_1\|_2\leq d\pi/2\sigma+r(\tau+1)$. Repeat these
constructions for $y$ to obtain an analogous vertex $y_1$. Then
  $$
  \|x-x_1\|_2,\:\|y-y_1\|_2 ~\leq~ \frac{d\pi}{2\sigma}+r(\tau+1),
  $$
and, by Theorem \ref{thm:diameter},
  $$
  d_G(x,x_1),\;d_G(y,y_1) ~\leq~ \left(
    \frac{d\pi}{2\sigma}+r(\tau+1)\right)\frac{K}{\lambda}.
  $$

By definition of $\tau$, there is a vertex $x_2$ that belongs to a
petal $R_{x_2}$ on $U_x$ (indeed, lying on the same side of $O$
along $U_x$) with $R_{x_2}\subseteq B_2^d(O,r(\tau+3/2))$. In the
worst case, all empty petals in $U_x$ occur non-consecutively
between $x_1$ and $x_2$, so for $n$ large enough, we have
  \begin{eqnarray*}
  d_G(x_1,x_2) & \leq &
     \frac{\|x_1\|_2}{r}+\frac{2Kr\tau}{\lambda} \\
  & \leq & \frac{\|x\|_2}{r}+\frac{d\pi}{2\sigma r}
      +\frac{2Kr\tau}{\lambda}.
  \end{eqnarray*}
The same construction goes through if we replace the $x$'s with $y$'s.
Moreover, $\|x_2-y_2\|_2 \leq (2\tau+3)r$,
so the shortest path in $G$ between $x_2$ and $y_2$ satisfies
  $$d_G(x_2,y_2) \leq 
    \frac{(2\tau+3)rK}{\lambda}.$$
Concatenating all the above paths, we find that
  \begin{eqnarray}
  d_G(x,y) & \leq & d_G(x,x_1)+d_G(x_1,x_2)+d_G(x_2,y_2)+
    d_G(y_2,y_1)+d_G(y_1,y) \nonumber \\
  & \leq & \frac{\|x\|_2+\|y\|_2}{r}+\left(
  \frac{d\pi}{\sigma}+(8\tau+5)r\right)\frac{K}{\lambda}
  +\frac{d\pi}{\sigma r}\nonumber \\
  & \leq & \frac{\diam_p(\mathbf{B})}{\diam_2(\mathbf{B})}
    \frac{\|x\|_2+\|y\|_2+O(\rho)}{\lambda}
    +O(\tau)+O\left(\frac{(\ln n)^{-1/d}}{\lambda}\right).
  \label{eqn:diamUBRough}
  \end{eqnarray}
By the definitions of $r$ and $\rho$ given in \eqref{eqn:defnR} and
\eqref{eqn:defnRho}, and the bounds on $\xi$ and $\tau$
(Lemmas~\ref{lem:volIntBound} and~\ref{lem:tauBound}),
it follows that
  $$
  \frac{\xi}{\vol(\mathbf{B})} ~\geq~ \lambda^d\frac{\ln \ln n}{\ln n}
    ~\geq~  \frac{\ln\ln n}{n},
  $$
where the second inequality follows from the assumption that $\lambda$ is
above the threshold for connectivity in Theorem~\ref{thm:conThresh}.  Therefore
  \begin{eqnarray*}
  \tau & \leq & \frac{2\sigma^{(d-1)/2}}{r}
  \exp{\left(\frac{-n\xi}{\vol(\mathbf{B})}\right)} \\
  & \leq & \frac{4}{\lambda}(\ln n)^{(d-1)/(2d)}(\ln n)^{-1} =
    \frac{O\left((\ln n)^{-(d+1)/(2d)}\right)}{\lambda}
    = \frac{o(\rho)}{\lambda}.
  \end{eqnarray*}
Plugging this bound for $\tau$ into (\ref{eqn:diamUBRough}) gives
  $$d_G(x,y)\leq \frac{\diam_p(\mathbf{B})}{\diam_2(\mathbf{B})}
    \frac{\|x\|_2+\|y\|_2+O(\rho)}{\lambda},
  $$
and the theorem follows.

We now explain the assumption \eqref{wolog-assumption}. If
$\|x\|_2,\|y\|_2\leq (\tau+3/2)r$, then by Theorem
\ref{thm:diameter}, $d_G(u,v)\leq (2\tau+3)rK/\lambda =
o(\rho)/\lambda$.  On the other hand, if $\|x\|_2\leq (\tau+3/2)r
\leq \|y\|_2$, then, almost always, there is a vertex $x_2\in
B_2^d(O,(\tau+3/2)r)$ belonging to some petal. By the preceding argument,
$d_G(x,y)\leq o(\rho/\lambda)+d_G(x_2,y)\leq
(\diam_p(\mathbf{B})/\diam_2(\mathbf{B}))(\|y\|_2+O(\rho))/\lambda.$
\end{proof}

We conclude with two remarks.  First, the proof of Theorem
\ref{thm:diameterBound} can easily be adapted to the case $d=1$ to
obtain the upper bound $\mathrm{diam}(G_p^1(\lambda,n)) \, \leq \,
\left(2+O(\ln \ln n/\ln n)\right)/\lambda$. (Note that when $d=1$
a pincushion consists of just one pin.) Second, the technique of
Theorem \ref{thm:diameterBound} can be extended to obtain the
stronger result $d_G(x,y)\leq \left(\|x-y\|_p+O\left((\ln \ln
n/\ln n)^{1/d}\right)\right)/\lambda$, so that the graph distance
\emph{approximates} the $\ell_p$-metric. Each pin is replaced by
$\sigma^{d-1}$ evenly spaced parallel pins. We can still bound
$\tau$ by $o(\rho)/\lambda$, but now any two vertices $x$ and $y$
are close to the same pin, on which a short path from $x$ to $y$
is found. We refer the reader to \cite{E05} for details.

\section*{Acknowledgements}

The authors thank Thomas Schlumprecht, Joel Spencer
and Dennis Stanton for various helpful discussions and comments.
We also thank the anonymous referees for several helpful
suggestions.

\appendix
\section{Facts about $\ell_p$- and spherical geometry}

In the body of the article, we used various facts about the
Euclidean and $\ell_p$-geometry of balls and spherical caps.
None of these facts are difficult; however, for convenience we
present them together here along with brief proofs.

\subsection{Volume of the $\ell_p$-unit ball}
\label{appendix:sphere-volume}

Fix $p\in[1,\infty]$ and an integer $d\geq 1$. Let $B^d_p(r)$ be
the $\ell_p$-ball centered at the origin in $\Rr^d$ with radius
$r$:
  $$B^d_p(r) = \left\{\vec{x} \in \Rr^d ~:~
    \sum_{i=1}^d |x_i|^p \leq r^p\right\}.$$
We will show that the ($d$-dimensional) volume of $B^d_p(r)$ is
  \begin{equation} \label{lp-ball-volume}
    \vol(B^d_p(r)) = \frac{(2r)^d \Gamma\left(\frac{p+1}{p}\right)^d}
         {\Gamma\left(\frac{p+d}{p}\right)},
  \end{equation}
where $\Gamma$ is the usual gamma function \cite{WhitWat}. For
$d=1$ this is trivial.  For $d>1$ we have
  $$
  \vol(B^d_p(r))
    = 2\int_0^r \vol\left(B^{d-1}_p((r^p-x^p)^{1/p})\right)~dx.
  $$
Make the substitution $u=x^p/r^p$,
$x=ru^{1/p}$, $dx=(r/p)u^{(1-p)/p}\,du$.
By induction on $d$, we obtain
  $$
  \vol(B^d_p(r))
    = \frac{2^dr^d\Gamma\left(\frac{p+1}{p}\right)^{d-1}}
             {p\Gamma\left(\frac{p+d-1}{p}\right)}
        \int_0^1 (1-u)^{\frac{d-1}{p}} u^{\frac{1-p}{p}}~du.
  $$
Evaluating this integral as in \cite[\S 12.4]{WhitWat} yields the
desired formula \eqref{lp-ball-volume}. It follows that
  \begin{equation} \label{eqn:a-d-p}
    \alpha_p^d ~:=~ \frac{\vol(B^d_p(r))}{\vol(B^d_2(r))}
   ~=~ \frac{\Gamma\left(\frac{p+1}{p}\right)^d \cdot \Gamma\left(\frac{2+d}{2}\right)}
            {\Gamma\left(\frac{3}{2}\right)^d   \cdot \Gamma\left(\frac{p+d}{p}\right)}\, .
  \end{equation}

\subsection{$\ell_p$-antipodes on the Euclidean sphere}
\label{appendix:antipodes}

Let $\mathbf{B}=B^d_2(1)$ be the Euclidean unit ball, centered at
the origin in $\Rr^d$, and let $p\neq 2$.  We wish to calculate
the $\ell_p$-diameter of $\mathbf{B}$, that is,
  $$\diam_p(\mathbf{B}) =
  \max\left\{\|\vec{x}-\vec{y}\|_p ~:~ \vec{x},\vec{y} \in
  \mathbf{B}\right\}.$$
A pair of points of $\mathbf{B}$ at distance $\diam_p(\mathbf{B})$
are called \emph{$\ell_p$-antipodes}.  If $d=1$, then $\mathbf{B}$
is a line segment and the only antipodes are its endpoints. Of
course, if $p=2$, then the antipodes are the pairs $\pm\vec{x}$
with $\|\vec{x}\|_2=1$.

Suppose that $d>1$.  Every pair of antipodes $\vec{x},\vec{y}$
must satisfy $\|\vec{x}\|_2=\|\vec{y}\|_2=1$. Without loss of
generality, we may assume $x_i\geq 0\geq y_i$ for every $i$. Using
the method of Lagrange multipliers (with objective function
$\left(\|\vec{x}-\vec{y}\|_p\right)^p = \sum_{i=1}^d
(x_i-y_i)^p$), we find that for every $i$,
  \begin{equation} \label{lagrange}
    p(x_i-y_i)^{p-1} = 2\lambda x_i = -2\mu y_i,
  \end{equation}
where $\lambda$ and $\mu$ are nonzero constants.  In particular
$y_i=(-\lambda/\mu)x_i$ for every $i$, so $\vec{y}=-\vec{x}$.
(That is, every pair of $\ell_p$-antipodes on $\mathbf{B}$ is a
pair of $\ell_2$-antipodes.) Substituting for $y_i$ in
\eqref{lagrange} gives
  $$\nu x_i^{p-1} = x_i,$$
where $\nu$ is some nonzero constant.  In particular, the set of
coordinates $\{x_1,\dots,x_d\}$ can contain at most one nonzero
value, so either $\vec{x}$ is a coordinate unit vector or else
$x_i=d^{-1/2}$ for every $i$.  In the first case
$\|(-\vec{x})-\vec{x}\|_p=2$, while in the second case
$\|(-\vec{x})-\vec{x}\|_p=2d^{-1/2+1/p}$. In summary:

\begin{proposition} \label{prop:antipodes}
Let $d\geq 2$, and let $\mathbf{B}=B^d_2(1)$ be the Euclidean unit
ball, centered at the origin in $\Rr^d$. Then the pairs of
$\ell_p$-antipodes on $\mathbf{B}$ are precisely the pairs
$\{\pm\vec{x}\}$ satisfying the following additional conditions:
  $$
  \begin{cases}
    |x_1| = \cdots = |x_d| = d^{-1/2}            &\quad\text{when } 1\leq p<2,\\
    \|\vec{x}\|_p = 1                              &\quad\text{when } p=2,\\
    \vec{x} \text{ is a coordinate unit vector}  &\quad\text{when } 2<p\leq \infty.
  \end{cases}
  $$
In particular, the $\ell_p$-diameter of $B$ is
  \begin{equation} \label{lp-diameter}
    \diam_p(B) ~=~ \max\left(2,\;2d^{1/p-1/2}\right) ~=~
    \begin{cases}
      2d^{1/p-1/2} & \quad\text{for } 1\leq p\leq 2,\\
      2 & \quad\text{for } 2\leq p\leq \infty.
    \end{cases}
  \end{equation}
\end{proposition}

\subsection{Spherical caps}
\label{appendix:spherical-cap-subsection}

Let $B=B^d_2(r)$ be the Euclidean ball of radius $r$, centered at
the origin in $\mathbb{R}^d$. Let $C$ be a spherical cap of $B$ of
height $h$, with $0\leq h\leq r$; for instance,
  $$C = \{(x_1,\dots,x_d)\in B:~ r-h\leq x_1\leq r\}.$$
The volume of $C$ can be determined exactly, but the precise
formula is awkward for large $d$ (one has to evaluate the integral
$\int\sin^d\theta~d\theta$). On the other hand, we can easily
obtain a lower bound for $\vol(C)$ by inscribing in it a
``hypercone'' $H$ of height $h$ whose base is a $(d-1)$-sphere of
radius $s=\sqrt{r^2-(r-h)^2}=\sqrt{2rh-h^2}$. For $r-h\leq x\leq
r$, the cross-section of $H$ at $x=x_1$ is $B_2^{d-1}(s(r-x)/h)$,
so applying \eqref{lp-ball-volume} gives
  \begin{eqnarray*}
    \vol(H)
    &=& \frac{2^{d-1}s^{d-1}\Gamma\left(\frac{3}{2}\right)^{d-1}}
             {h^{d-1}\Gamma\left(\frac{d+1}{2}\right)}
        \int_{r-h}^r (r-x)^{d-1}~dx\\
    &=& \frac{\pi(2r-h)^{(d-1)/2} h^{(d+1)/2}}{d\cdot\Gamma\left(\frac{d+1}{2}\right)}.
  \end{eqnarray*}
Since $2r-h\geq r$, we have the bound
  \begin{equation} \label{cap-bound}
    \vol(C) \geq \frac{\pi r^{(d-1)/2}}{d\cdot\Gamma\left(\frac{d+1}{2}\right)} \cdot h^{(d+1)/2}.
  \end{equation}

\subsection{The $\ell_p$-distance between opposite spherical caps}
\label{appendix:cap-distance}

\begin{figure}
\includefigure{4in}{2in}{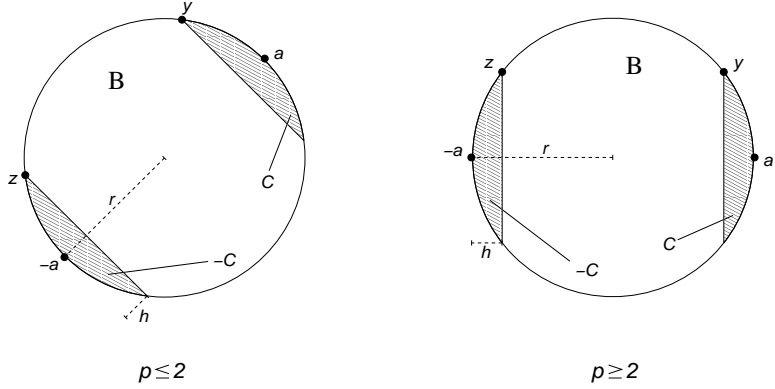}
\caption{Antipodally aligned caps $\pm C$ centered at antipodes
$\pm a$ of $B\subseteq \mathbb{R}^2$ for different values of $p$
below and above $2$.\label{fig:antipodes}}
\end{figure}

By definition, the $\ell_p$-distance between two sets
$Y,Z\subseteq \Rr^d$ is
  $$d_p(Y,Z):=\inf\{\|y-z\|_p~:~y\in Y,~z\in Z\}.$$
Let $B=B^d_2(r)$, let $\pm a$ be a pair of $\ell_p$-antipodes on
$B$, and let $\pm C$ be the cap of height $h$ centered at $\pm a$.
Note that every pair $y,z$ at minimum distance has displacement
parallel to $a$ (see Figure~\ref{fig:antipodes}).  (This can be
verified by another easy Lagrange-multiplier calculation.)
Therefore
  \begin{equation} \label{cap-distance-formula}
    d_p(C,-C) ~=~ \frac{2(r-h)}{r}\|a\|_p ~=~ \begin{cases}
      2(r-h)d^{1/p-1/2} & \quad\text{for } 1\leq p\leq 2,\\
      2(r-h) & \quad\text{for } p\geq 2.
    \end{cases}
  \end{equation}

\end{document}